\documentclass[11pt]{article}
\usepackage{amsmath, amsfonts, amssymb}
\usepackage{graphicx, amsthm,color}

\providecommand{\U}[1]{\protect\rule{.1in}{.1in}}

\newtheorem{theorem}{Theorem}[section]
\newtheorem{lemma}[theorem]{Lemma}

\newtheorem{corollary}[theorem]{Corollary}
\newtheorem{proposition}[theorem]{Proposition}
\theoremstyle{definition}

\theoremstyle{remark}
\newtheorem{remark}[theorem]{Remark}
\numberwithin{equation}{section}

\def \R{{\bf R}}

\def \dist{\textrm{dist}}

\def \R{\mathbb{R}}

\begin{document}
	
	\begin{center}
		{\LARGE The slope robustly determines convex functions}
		
		\bigskip
		
		\textsc{Aris Daniilidis \& Dmitriy Drusvyatskiy}
	\end{center}
	
	\bigskip
	
	\noindent\textbf{Abstract.} 
	We show that the deviation between the slopes of two convex functions controls the  deviation between the functions themselves.	This  result reveals that the slope---a one dimensional construct---robustly determines convex functions, up to a constant of integration. 
	\bigskip
	
	
	\vspace{0.6cm}
	
	\textit{Secondary} 37C10, 49J52.
	
	\section{Introduction}
The recent paper~\cite[Theorem~3.8]{BCD2018} established the following intriguing result. Two
$\mathcal{C}^{2}$-smooth, convex and bounded from below functions $f,g$
defined on a Hilbert space $\mathcal{H}$ are equal up to an additive constant
if and only if their gradient norms coincide:
\begin{equation}
	\|\nabla f\|=\|\nabla g\|\qquad\;\Longleftrightarrow\;\qquad f=g+\mathrm{cst.}\label{eq:det}
\end{equation}
This result is ostensibly surprising since it readily yields that the function $x\mapsto \|\nabla f(x)\|$, which takes values in the real line, determines the entire gradient map $x\mapsto \nabla f(x)$, which takes values in~$\mathcal{H}$. In the follow up work \cite{PSV2021}, the assumption on smoothness of $f$ was further weakened to continuity with the gradient norm $\|\nabla f(x)\|$ replaced by the slope $s_f(x):=\dist(0,\partial f(x))$. Here $\partial f(x)$ denotes the subdifferential of the convex function $f$ at $x$.\footnote{We note that further
	generalizations of the determination result \cite{PSV2021}
	have recently been achieved: for convex continuous bounded from below
	functions in Banach spaces (see~\cite{TZ2022}) and for Lipschitz coervice
	functions in metric spaces (\cite{DS2022}). For the time being, we do not
	pursue our sensitivity analysis in this generality.}

In this work, we ask whether the slope (or the gradient norm in the smooth case) {\em robustly} determines the function itself. That is, if the slopes for two functions are close, then how close are the function values? Roughly speaking, we will show that for any two continuous convex functions $f$ and $g$ defined on a Hilbert space, the following estimate is true:
$$\|g-f\|_{\mathcal{U}}\,\,\lesssim \,\, \|s_g-s_f\|_{\mathcal{U}}\,+\,\sqrt{\|s_g-s_f\|_{\mathcal{U}}} \,+\,\|g-f\|_{C_f\,\cup\, C_g}.$$
Here $\mathcal{U}$ is any bounded set where $f$ is bounded, $\|\cdot\|_{\mathcal{U}}$ denotes the sup-norm over $\mathcal{U}$, and $C_f$ and $C_g$ are the sets of minimizers of $f$ and $g$, respectively. In particular, the deviation $\|g-f\|_{\mathcal{U}}$ exhibits a dependence on $\|s_g-s_f\|_{\mathcal{U}}$ that is at worst H\"{o}lder with exponent $1/2$. In the finite-dimensional setting $\mathcal{H}=\R^n$, we show that this undesirable square root dependence  may be dropped:   
$$\|g-f\|_{\mathcal{U}}\lesssim  \|s_g-s_f\|_{\mathcal{U}}+\|g-f\|_{C_f\cup C_g}.$$
The downside is that the hidden constant in this bound depends on the length of subgradient curves initialized in $\mathcal{U}$ and at worst grows super exponentially in the dimension $n$.

	\section{Notation and preliminaries}
	
	Let $\mathcal{H}$ denote a Hilbert space and let $f:\mathcal{H}\rightarrow
	\mathbb{R}$  be  a convex continuous function. We denote the set of minimizers of $f$ by
	\[
	\mathcal{C}_{f}:=\arg\min f,
	\]
	and suppose that $\mathcal{C}_{f}$ is nonempty (therefore the infimum value $f_{\ast}:=\inf f$ is attained). The key object we will focus on is the slope $s_f(x)=\dist(0,\partial f(x)),$ where $\partial f(x)$ denotes the subdifferential: 	
\begin{equation}\label{eq:fenchel}
\partial f(x)=\{v\in \mathcal{H}:\, f(y)-f(x)\geq \langle v,y-x\rangle,~~\forall x,y\in \mathcal{H}\}.
\end{equation}
Equivalently, $s_f(x)$ measures the fastest instantaneous rate of decrease of $f$ from $x$. 

Our goal is to show that the deviation between the slopes of two convex functions controls the deviation between the functions themselves.	Our arguments will make heavy use of subgradient dynamical systems, a topic we review now following~\cite{ABM2014-book,Brezis}. Namely, \cite[Theorem~17.2.2]{ABM2014-book} shows that for every initial point $x\in\mathcal{H}$, there exists a unique, maximally defined, injective, absolutely continuous curve $\gamma:[0,T_{\max})\rightarrow
	\mathcal{H}$, such that
	\begin{equation}\label{eqn:subgrad_system}
		\left\{
		\begin{array}
			[c]{l}
			\dot{\gamma}(t)\underset{\text{a.e}}{\in}-\partial f(\gamma(t))\medskip\\
			\gamma(0)=x
		\end{array}
		\right.  \tag{GS}
	\end{equation}
Subgradient curves $\gamma$ satisfy a number of useful properties, summarized below.
 \begin{enumerate}
\item[(P1)] Equality 
	\begin{equation}
		\|\dot{\gamma}(t)\|=s_{f}(\gamma(t))  \quad\text{holds for {a.e.}
		}t\in\lbrack0,T_{\max}). \label{eq:s1}
	\end{equation}
and the slope function $t\mapsto s_{f}(\gamma(t))$ is nonincreasing on $[0,T_{\max})$.

\item[(P2)] The function $r(t)=f(\gamma(t))$ is
		convex and strictly decreasing on  $[0,T_{\max})$, and
		\[
		\underset{t\rightarrow T_{\max}}{\lim}f(\gamma(t))=f_{\ast}.
		\]
\item[(P3)] The distance function $t\mapsto d(\gamma(t),\mathcal{C}_{f})$ is
		strictly decreasing on $[0,T_{\max})$. Moreover, for
		every $x_{\ast}\in\mathcal{C}_{f},$ the function $t\mapsto\Vert\gamma(t)-x_{\ast}\Vert$ is
		strictly decreasing on $[0,T_{\max})$.
\end{enumerate}	
Property (P1) follows from~\cite[Theorem~17.2.2~(iii)-(iv)]{ABM2014-book}, (P2) is given in~\cite[Proposition~17.2.7~(i)]{ABM2014-book}, while (P3) follows easily after differentiation, using (GS) and~\eqref{eq:fenchel}.  \smallskip\newline
Next, we will require two estimates on the length of subgradient curves. The first (Lemma~\ref{Lemma_length}) is an easy consequence of (P1) and (P2) above (we provide a proof for convenience), while the second (Proposition~\ref{Prop_SCC_length}) was essentially proved in~\cite{MP1991} for a particular class of Lipschitz curves (therein called $\Gamma$-curves, ultimately known as \textit{self-contracted} curves, definition coined in~\cite{DLS2010}) and became explicit for subgradient curves in \cite{DDDL2015,LMV2015}.

	\begin{lemma}
		[Length estimation I]\label{Lemma_length} Let $f\colon
		\mathcal{H}\rightarrow\mathbb{R}$ be a convex continuous function with
		nonempty set of minimizers and let
		$\gamma:[0,T_{\max})\rightarrow\mathcal{H}$ be the solution of (GS). Then for
		every $T\in(0,T_{\max})$, setting $\gamma_{T}:=\gamma(T)$ we have:
		\[
		\int_{0}^{T}
		\,|\dot{\gamma}(t)|\,dt\,\leq\left[  s_{f}(\gamma_{T})\right]  ^{-1}\left(
		f(x)-f_{\ast}\right)  .
		\]
		
	\end{lemma}
	
	\noindent\textbf{Proof. } Set $r(t):=f(\gamma(t))$ and denote by
	$h\ $the inverse function of the mapping $t\mapsto r(t)$ on the interval $[0,T_{\max})$. Then for the
	reparametrization $\tilde{\gamma}(\rho)=\gamma(h(\rho))$ we have $f(\tilde{\gamma}(\rho))=\rho$. Differentiating gives
	\[
	\frac{d}{d\rho}[\tilde{\gamma}(\rho)]=\frac{\partial f(\tilde\gamma(\rho))^{\circ}}{s_{f}(\tilde{\gamma}(\rho))^{2}},\quad\text{for \textit{a.e. }}
	\rho\in(f_*,f(x)],
	\]
	where $\partial f(\tilde\gamma(\rho))^{\circ}$ is the element of $\partial f(\tilde\gamma(\rho))$ of minimal norm, thus 
$\|\partial f(\tilde\gamma(\rho))^{\circ}\|=s_f(\tilde\gamma(\rho))$.
	Taking into account that the function $\rho\mapsto s_{f}(\tilde{\gamma}
	(\rho))$ is increasing, we deduce:
	\[
	{\displaystyle\int\limits_{0}^{T}}
	\,\|\dot{\gamma}(t)\|\,dt\,=
	{\displaystyle\int\limits_{f(\gamma_{T})}^{f(x)}}
	\,\frac{1}{s_{f}(\tilde{\gamma}(\rho))}\,d\rho\,\,\leq\,\,\frac{f(x)-f(\gamma_{T}
		)}{s_{f}(\gamma_{T})}
	\]
	and the result follows.\hfill$\Box$
	
	\bigskip
	
\begin{proposition}[Length estimation II] \label{Prop_SCC_length} Assume
$\mathcal{H}=\mathbb{R}^{n}$. There exists a constant $K_{n}$ depending only on dimension such that for every $x\in\mathbb{R}^{n}$ the solution $\gamma(\cdot)$ of the subgradient system~\eqref{eqn:subgrad_system}
has length bounded by $K_{n}\cdot d(x,\mathcal{C}_{f})$.
\end{proposition}

\noindent The above result provides a universal bound $K_n$ for the ratio between the length of a subgradient curve and its diameter, the drawback being that that the dependence of $K_n$ on the dimension is of the order of $n^{n/2+1}$ (see \cite{MP1991,GBR2021}).	
	
\section{Main results}

	For any function $\omega\colon\mathcal{H}\to \R$ and a set $\mathcal{U}\subset\mathcal{H}$, we will use the notation
	$$\Vert \omega|_{\mathcal{U}}:=\,\underset{x\in
			{\mathcal{U}}}{\sup}~ \left(\max\,\{\omega(x),0\}\right)\qquad \textrm{and}\qquad \|\omega\|_{\mathcal{U}}=\underset{x\in {\mathcal{U}}}{\sup}~ |\omega(x)|.$$
Note that $\Vert \omega|_{\mathcal{U}}$ provides a one-sided bound\footnote{Notice that $\Vert \cdot|_{\mathcal{U}}$ is the canonical asymmetrization of the seminorm $\Vert \cdot\Vert_{\mathcal{U}}$ of uniform convergence, see \cite{DSV2021}.}, while  $\|\omega\|_{\mathcal{U}}$ is the standard two-sided sup-norm. \smallskip\newline The following is the main theorem of the paper.
	
	\begin{theorem}\label{Theorem_main}Let $f,g\colon \mathcal{H}\rightarrow\mathbb{R}$  be convex continuous functions. Assume $\mathcal{C}_{f}=\arg\min f\neq\emptyset$ and set $f_{\ast}=\min f$.	
			 For each $r>0$ define the tube around $C_f$ by 
		\begin{equation}
			\mathcal{U}_{r}:=\{x\in\mathcal{H}:\,d(x,\mathcal{C}_{f})\leq r\}.
		\end{equation}
Then for every $x\in \mathcal{U}_{r}$, the estimate holds:
\begin{equation}\label{eqn:main_ineq}
g(x)-f(x)\,\leq \,\Vert s_{g}-s_{f}|_{\mathcal{U}_{r}}\,+\,\Vert g-f|_{\mathcal{C}_{f}}\,+\,2\sqrt{d(x,\mathcal{C}_{f})\cdot\Vert s_{g}-s_{f}|_{\mathcal{U}_{r}}\cdot \left(  f(x)-f_{\ast}\right)}.
\end{equation}
Moreover, in the finite-dimensional setting $\mathcal{H}=\R^n$, there exists a constant $K_n>0$ depending only on the dimension $n$ such that 
\begin{equation}
	g(x)-f(x)\,\leq \, K_{n}\,\,\Vert s_{g}-s_{f}|_{\mathcal{U}_{r}}\,\,d(x,\mathcal{C}_{f})\,+\,\Vert
	g-f|_{\mathcal{C}_{f}}. \label{eq:cv1}
\end{equation}
	\end{theorem}
	
	\noindent\textbf{Proof. }Let $x\in\mathcal{H}\setminus\mathcal{C}_f$ be arbitrary and fix $\delta>0$. 	Our goal is to show the estimate
	\begin{equation}
		g(x)-f(x)\,\leq\left(\Vert s_{g}-s_{f}|_{\mathcal{U}_r}+\delta\right)\,d(x,\mathcal{C}_{f})\,+\,\frac{\Vert s_{g}-s_{f}|_{\mathcal{U}_r}}{\delta}\,\,\left(  f(x)-f_{\ast}\right)\,+\,\Vert g-f|_{\mathcal{C}_{f}}, \label{eq:ad1}
	\end{equation}
	from which \eqref{eqn:main_ineq} follows by setting
	$\delta=\sqrt{\frac{\Vert s_{g}-s_{f}|_{\mathcal{U}_r}\cdot \left(  f(x)-f_{\ast}\right)}{d(x,\mathcal{C}_{f})}}.$\smallskip\newline
We consider two cases:\smallskip

(i). Suppose that $s_f(x)\leq \delta$ and let $\hat{x}:=\mathrm{proj}_{\mathcal{C}_{f}}(x)$ be the projection of $\hat{x}$ to the closed convex set~$\mathcal{C}_{f}$ (therefore $f(\hat{x})=f_*\leq f(x)$). Then we compute
	\[
g(x)-g(\hat{x})\leq s_{g}
	(x)\,\|x-\hat{x}\|\,\leq(\Vert s_{g}-s_{f}|_{\mathcal{U}_r}+\delta)\,d(x,\mathcal{C}_{f}),
	\]
	where the first inequality follows from convexity of $g$.
	We therefore conclude
	\begin{align*}
	g(x)-f(x)&=(g(x)-g(\hat x))+(g(\hat x)-f(\hat x))+(f(\hat x)- f(x))\\
	&\leq (\Vert s_{g}-s_{f}|_{\mathcal{U}_r}+\delta)\,d(x,\mathcal{C}_{f}) +\Vert g-f|_{\mathcal{C}_{f}},	\end{align*}
	thus verifying \eqref{eq:ad1}. \smallskip

(ii). Suppose now that $s_f(x)> \delta$ and let $\gamma\colon \lbrack0,T_{\max}) \to \mathcal{H}$ denote the unique maximal solution of the subgradient system \eqref{eqn:subgrad_system} for $f$.
Define the function
\[
a(t):=f(\gamma(t))-g(\gamma(t)).
\]
Differentiating, for \textit{a.e.} $t\in\lbrack0,T_{\max})$, we have (\textit{c.f.} \cite[Proposition~17.2.5]{ABM2014-book}):
\[
\dot{a}(t)=-s_{f}(\gamma(t))^{2}-\langle\partial g(\gamma(t))^{\circ},\dot{\gamma}(t)\rangle,
\]
where $\partial g(\gamma(t))^{\circ}$ is the element of minimal norm of  $\partial g(\gamma(t))$, that is, $s_{g}(\gamma(t))=\Vert\partial g(\gamma(t))^{\circ}\Vert$. From the Cauchy-Schwarz inequality we conclude:
\begin{equation}\label{eqn:basic}
	\begin{aligned}
	\dot{a}(t)&\leq -s_{f}
	(\gamma(t))^{2}+s_{g}(\gamma(t))\cdot s_{f}(\gamma
	(t))\\
	&= \left(  s_{g}(\gamma(t))-s_{f}(\gamma(t))\right)  \,s_{f}(\gamma
	(t))\,\\
	&\leq\Vert s_{g}-s_{f}|_{\mathcal{U}_r}\,\|\dot{\gamma}(t)\|.
	\end{aligned}
\end{equation}
Define
	\[
	T:=\sup\ \{t\in\lbrack0,T_{\max}):\,s_f(\gamma(t))>\delta\}.
	\]
Setting
	$\gamma_{T}:=\gamma(T)$ and integrating~\eqref{eqn:basic} on~$[0,T]$ we obtain:
\begin{equation}
g(x)\,\leq\, f(x)\,+\,\left[  g(\gamma_{T})-f(\gamma_{T})\right] \, +\,\Vert s_{g}-s_{f}|_{\mathcal{U}_r}\,\int_{0}^{T}\!\!\|\dot{\gamma}(t)\|\,dt. \label{eq:ad3}
\end{equation}
By Lemma~\ref{Lemma_length} and the definition of $T$ we get:
	\begin{equation}
		\int_{0}^{T}\!\|\dot{\gamma}(t)\|\,dt\,\leq\,\left[  s_{f}(\gamma
		_{T})\right]  ^{-1}\left(  f(x)-f_{\ast}\right)  \leq\delta^{-1}\,\left(
		f(x)-f_{\ast}\right)  . \label{eq:ad4}
	\end{equation}
Let $\hat{\gamma}=\mathrm{proj}_{\mathcal{C}_{f}}(\gamma_{T})$ be the projection of $\gamma_{T}$ to the set of minimizers 
$\mathcal{C}_{f}$. Then
	\[
	f(\hat{\gamma})=f_{\ast}\leq f(\gamma_{T})\qquad\text{and}\qquad\Vert
	\gamma_{T}-\hat{\gamma}\Vert=d(\gamma_{T},\mathcal{C}_{f})\leq d(x,\mathcal{C}_{f}).
	\]
Taking into account $s_f(\gamma_{T})\leq \delta$ we deduce $s_{g}(\gamma_{T})\leq\Vert s_{g}-s_{f}|_{\mathcal{U}_r}+\delta$ and consequently
	\[
	g(\gamma_{T})-g(\hat{\gamma})\leq s_{g}(\gamma_{T})\,\|\gamma_{T}-\hat{\gamma}\|\,\leq\left(\Vert s_{g}-s_{f}|_{\mathcal{U}_r}+\delta\right)\,d(x,\mathcal{C}
	_{f}),
	\]
	where the first inequality follows from convexity of $g$.
	We readily obtain that:
	\begin{equation}
		g(\gamma_{T})-f(\gamma_{T})\,\leq\,\left(  g(\gamma_{T})-g(\hat{\gamma})\right)
		+\left(  g(\hat{\gamma})-f_{\ast}\right) \, \leq\left(\Vert s_{g}-s_{f}|_{\mathcal{U}_r
		}+\delta\right)d(x,\mathcal{C}_{f})\,+\,\Vert g-f|_{\mathcal{C}_{f}}.
		\label{eq:ad5}
	\end{equation}
	Combining \eqref{eq:ad3}, \eqref{eq:ad4}, and \eqref{eq:ad5} yields the claimed estimate~\eqref{eq:ad1}.
	Finally, the estimate \eqref{eq:cv1} follows by letting $T\uparrow T_{\max}$ in \eqref{eq:ad3} and using Proposition~\ref{Prop_SCC_length} to bound the length of $\gamma(\cdot)$. 
 \hfill$\Box$
	
\bigskip
	
\noindent An easy consequence of the above is the following guarantee of asymptotic consistency.
	
	\begin{corollary}[Robust (one-sided) determination]
		\label{Cor_Lip_seq} Let $f, \{f_{k}\}_{k\geq 0}:\mathcal{H}\rightarrow\mathbb{R}$
		be convex continuous functions and suppose that $\mathcal{C}_{f}$ is nonempty and bounded. 
		Assume further that
\begin{itemize}
\item[{\em (i).}] $\underset{k\geq1}{\lim\sup }\,\|s_{f_{k}}- s_{f}|_{\mathcal{U}}\leq 0$, for all  bounded sets $\mathcal{U}\subset\mathcal H$; and			
			\item[{\em (ii).}] $\underset{k\geq1}{\lim\sup }\,\|f_{k}- f |_{C_f}\leq 0$.
\end{itemize}
		
		\noindent Then $\underset{k\geq1}{\lim\sup }\,\| f_{k}-f|_{\mathcal{U}} \leq 0$ for all bounded sets 
$\mathcal{U}\subset\mathcal{H}$.
	\end{corollary}
	
	\noindent\textbf{Proof.} Recalling from~Theorem~\ref{Theorem_main} the definition of $\mathcal{U}_{r}$, we
	observe that $\mathcal{U}_{r}$ is bounded. Our assumption can then be restated
	as follows:
	\[
	\forall r>0:\;\underset{k\geq1}{\lim\sup }\,\Vert s_{f_{k}}-s_{f}|_{\mathcal{U}_{r}}\leq0\qquad\text{and}\qquad\underset{k\geq 1}{\lim\sup }\,\Vert f_{k}-f|_{\mathcal{C}_{f}}\leq0.
	\]
	An application of Theorem~\ref{Theorem_main} for each $r>0$ completes the proof. \hfill$\Box$

\bigskip

\noindent A symmetric version of the corollary follows by an analogous argument. 
\begin{corollary}
		[Robust (two-sided) determination]\label{Cor_Lip_sym}Let
		$f,\{f_{k}\}_{k\geq1}\colon\mathcal{H}\rightarrow\mathbb{R}$ be convex
		continuous functions such that
		\[
		\mathcal{C}_{f_{k}}\neq\emptyset,\;\forall k\geq1\qquad\text{and\qquad
		}\mathcal{C}:=\mathcal{C}_{f}\cup\left(  \cup_{k\geq1}\mathcal{C}_{f_{k}}\right)  \text{\ is bounded.}
		\]
		Assume further that:\smallskip

{\em (i).} $s_{f_{k}}$ converge to $s_{f}$ uniformly on bounded sets,
		  \smallskip

{\em (ii).} $f_{k}$ converge to $f$ uniformly on $\mathcal{C}$.\medskip\newline 
		Then $f_{k}$ converge to $f$ uniformly on bounded sets.
	\end{corollary}

\begin{remark}[open question] Our approach is heavily based on the existence of minimizers. We do not know if the results of this work can be extended to the class of lower semicontinuous convex functions, which are bounded for below. This is a challenging question that merits investigation. 
\end{remark}
\medskip

\begin{center}
\noindent\rule{4cm}{2pt}
\end{center}

	\medskip
	
	\noindent\rule{5cm}{1pt} \smallskip\newline
	\noindent\textbf{Acknowledgements.} A major part of this work has been
	accomplished during a research visit of the first author to the University of
	Washington. This author thanks the host institution for hospitality and
	acknowledges support from the Austrian Science Fund (FWF, P-36344-N).

\newpage
	
	\noindent Aris DANIILIDIS
	
	\medskip
	
	\noindent Institute of Statistics and Mathematical Methods in Economics,
	VADOR E105-04 \newline TU Wien, Wiedner Hauptstra{\ss }e 8, A-1040 Wien\medskip
	\newline(on leave) DIM-CMM, CNRS IRL 2807 \newline Beauchef 851, FCFM,
	Universidad de Chile 
           \medskip\newline\noindent E-mail:
	\texttt{aris.daniilidis@tuwien.ac.at}\newline\noindent
	\texttt{https://www.arisdaniilidis.at/}
	
	\medskip
	
	\noindent Research supported by the grants: \smallskip\newline Austrian Science Fund (FWF P-36344N) (Austria)\\ 
	CMM  FB210005 BASAL funds for centers of excellence (ANID-Chile)\newline
	
	\bigskip
	
	\noindent Dmitriy Drusvyatskiy
           \medskip
           \newline University of Washington \newline
	Department of Mathematics \newline C-138 Padelford, Seattle, WA 98195
           \medskip\newline
	\noindent E-mail: \texttt{ddrusv@uw.edu} \newline\noindent
	\texttt{http://www.math.washington.edu/{\raise.17ex\hbox{$\scriptstyle\sim$}}ddrusv/}. \smallskip

\end{document}